\documentclass[11pt]{article}
\usepackage{latexsym,euscript,amsmath,amssymb,amsbsy,amsfonts,amsthm,amsopn,amstext,amsxtra,amscd}
\usepackage{epsfig}
\usepackage{graphics}
\usepackage{url}

\usepackage[english]{babel}

\theoremstyle{plain}
\newtheorem{theorem}{Theorem}

\newtheorem{corollary}{Corollary}

\theoremstyle{definition}
\newtheorem{example}{Example}

\newtheorem{remark}{Remark}

\newcommand{\N}{{\mathbb N}}

\newcommand{\Z}{{\mathbb Z}}

\newcommand{\D}{\text{$\mathcal{D}$}}

\begin{document}
\centerline{\bf ON DISTRIBUTION OF CONTINUOUS SEQUENCES}
\centerline{Milan Pa\v st\'eka} \vskip1cm Let $\N$ be the set of
positive integers and $\Omega$ compact metric ring of polyadic
integers.(see [N], [N1]) Let us remark that $\Omega$ is the
completion of $\N$ with respect to polyadic metric
\begin{equation}
\label{metric} \mathfrak{d}(a,b) =\sum_{n=1}^\infty
\frac{\psi_n(a-b)}{2^n},
\end{equation}
where $\psi(x)=0$ if $n|x$ and $\psi(x)=1$ otherwise. We shall use
two synonymous : sequence of real numbers and arithmetic function. A
sequence of real numbers $\{v(n)\}$ is called {\it polyadicly
continuous} if and only for each $\varepsilon >0$ such $m \in \N$
exists that
$$
\forall a, b \in \N; a\equiv b \pmod{m} \Rightarrow |v(a)-v(b)| <
\varepsilon.
$$
It is well known that polyadicly continuous sequence of real numbers
is uniformly continuous with respect polyadic metric $\mathfrak{d}$
and so each  polyadicly continuous sequence of real numbers
$\{v(n)\}$ can be extended by the natural way to a real valued
continuous function $\tilde{v}$ defined on $\Omega$ such that
$$
\tilde{v}(\alpha) = \lim_{j\to \infty}v(n_j)
$$
where $\{n_j\}$ is such sequence of positive integers that $n_j \to
\alpha$ for $j \to \infty$ with respect the polyadic metric. The
compact ring $\Omega$ is equipped with Haar probability measure $P$
and so the function $\tilde{v}$ can be considered as random variable
in the probability space $(\Omega, P)$. As usually if $h$ ia a
random variable on $\Omega$ the we denote $E(h)=\int h dP$ the mean
value of $h$.

Let $m \in \N$ and $s = 0,1, \dots m-1$. Put
$$
s+m\Omega = \{s+m\alpha; \alpha \in \Omega\}
$$
The ring $\Omega$ can be represented as disjoint decomposition
$$
\Omega = \bigcup_{s=0}^{m-1} s +m\Omega.
$$
(see [N], [N1]).  And so for the Haar probability measure $P$ we
have
\begin{equation}
\label{mesure} P(s +m\Omega)= \frac{1}{m}
\end{equation}
for $m \in \N, s=0, \dots m-1$.

We say that a set $S \subset \N$ has asymptotic density if and only
if the limit
$$
\lim_{n \to \infty} \frac{1}{n} card (S \cap [1,n]) := d(S)
$$
exists. The value $d(S)$ is then called {\it asymptotic density} of
$S$. {see [NAR], [G], [G1]). The system of all subsets of $\N$
having asymptotic density we shall denote $\D$. The function $d$ is
a finitely additive probability measure on $\D$. If $r+(m)$, $m \in
\N, r \in \Z$ is the set of all positive integers congruent to $r$
modulo $m$, then it is easy to check that $r+(m) \in \D$ and
$$
d(r+(m))=\frac{1}{m}.
$$

By the symbol $\mu^\ast(S)$ we shall denote the {\it Buck's measure
density } of the set $S \subset \N$ constructed in [BUC] as follows
$$
\mu^\ast(S) = \inf \Big\{\sum_{\ell=1}^k \frac{1}{m_\ell}; S \subset
\bigcup_{\ell=1}^k r_\ell +(m_\ell) \Big\}
$$
for $S \subset \N$.

In [PAS4] is proved
\begin{equation}
\label{clo}
 \mu^\ast(S) = P(cl(S))
\end{equation}
for each $S \subset \N$, where $cl(S)$ denote the topological
closure of $S$ in $\Omega$. It is known that $\mu^\ast$ is a strong
submeasure , the system of Buck's measurables sets
$$
\D_\mu=\{S \subset \N; \mu^\ast(S)+\mu^\ast(\N \setminus S) \}
$$
is an algebra of sets and the restriction $\mu=\mu^\ast|_{\D_\mu}$
is a finitely additive probability measure on $\D_\mu$. Moreover
$\D_\mu \subset \D$ and $\mu(S)=d(S)$ for each $S \in \D_\mu$.

The notion of uniform distribution of sequences was introduced and
studied first time by H. Weyl in his work [Wey]. $\{n \in \N; v(n)
<x\}$ belongs to $\D$ and
$$
d(\{n \in \N; v(n) <x\})=x.
$$

In the paper [PAS2], (see also [PAS3]) is this concept transfer for
the case of Buck's measure density and Buck's measurability.

 We say that a sequence $\{v(n)\}$ is {\it Buck's
measurable} if and only if for every $x$ - real number the set
$\{n\in \N; v(n) < x\}$ belong to $\D_\mu$.

A Buck's measurable sequence $\{v(n)\}$ is called {\it Buck's
uniformly distributed} if and only if for the function
$$
F(x)= \mu(\{n\in \N; v(n) < x\})
$$
we have
\begin{equation}\label{bud}
F(x)=0, x<0, F(x)=x, x \in [0,1] , F(x)=1, x>1.
\end{equation}

\begin{example}
\label{vdc} Let $\{Q_k\}$ be an increasing sequence of integers such
that $Q_0=1$ and $Q_k | Q_{k+1}, k=1,2,3,\dots$. Each positive
integer $n$ can be uniquely represented in the form
$$
n = a_0+ a_1Q_1+ ...+ a_kQ_k
$$
where $a_j < \frac{Q_{j+1}}{Q_j}, j=1,\dots ,k$. To this $n$ can be
associated an element of unit interval $\gamma(n)$ in the form
$$
\gamma(n)= \frac{a_0}{Q_1}+\dots + \frac{a_k}{Q_{k+1}}.
$$
The sequence $\{\gamma(n)\}$ is known as van der Corput sequence
with base $\{Q_k\}$ and in [PAS2] is proved that it is Buck's
uniformly distributed and polyadicly continuous.
\end{example}

\begin{theorem}
\label{thm1} {\bf Let $\{v(n)\}$ a polyadicly continuous sequence of
real numbers. Suppose that $F: (-\infty,  \infty) \to [0,1]$ is a
continuous function. If $F$ is the distribution function of random
variable $\tilde{v}$ then for each $x$ - real number we have that
$\{n\in \N; v(n) < x\} \in \D_\mu$ and}
\begin{equation}
\label{e1} \mu(\{n\in \N; v(n) < x\}) = F(x).
\end{equation}
\end{theorem}

{\bf Proof.} From the continuity of $F$ we get
\begin{equation}
\label{e2} P(\tilde{v} < x) = F(x)=P(\tilde{v} \le x)
\end{equation}
for each $x$ - real number. From the inclusion
$$
\{n\in \N; v(n) < x\} \subset \{\alpha \in \Omega; \tilde{v}(\alpha)
\le x\}
$$
we obtain
$$
cl(\{n\in \N; v(n) < x\}) \subset \{\alpha \in \Omega;
\tilde{v}(\alpha) \le x\}.
$$
And (\ref{e2}) yields
$$
\mu^\ast(\{n\in \N; v(n) < x\}) \le F(x)
$$
for every real number $x$. From the other side
$$
\N \setminus \{n\in \N; v(n) < x\} = \{n\in \N; v(n) \ge x\},
$$
therefore
$$
cl(\N \setminus \{n\in \N; v(n) < x\}) \subset\{\alpha \in \Omega;
\tilde{v}(\alpha) \ge x\}.
$$
Thus
$$
\mu^\ast(\N \setminus \{n\in \N; v(n) < x\}) \le 1-F(x)
$$
and the assertion follows. \qed

\begin{theorem}
\label{thm2} {\bf Let $\{v(n)\}$ be a polyadicly continuous sequence
and $F: (-\infty,\infty) \to [0,1]$ a continuous function. If for
every $x$ - real number
\begin{equation}
\label{e3} \mu^\ast(\{n\in \N; v(n) < x\})=F(x)
\end{equation}
holds, then the random variable $\tilde{v}$ has the distribution
function $F$ and for $x$ - real number the set $\{n\in \N; v(n) <
x\}$ belongs to $\D_\mu$ and}
$$
\mu(\{n\in \N; v(n) < x\})=F(x).
$$
\end{theorem}

{\bf Proof.} Clearly
$$
\{n\in \N; v(n) < x\} \subset \{\alpha \in \Omega; \tilde{v}(\alpha)
\le x\})
$$
and so $F(x) \le P(\tilde{v} \le x)$. From the other side
$$
\{\alpha \in \Omega; \tilde{v}(\alpha) < x\}) \subset cl(\{n\in \N;
v(n) \le x\} )
$$
thus for each $\varepsilon>0$ we have $F(x) \le P(\tilde{v} \le x)
\le F(x+\varepsilon)$. For $\varepsilon \to 0^+$ we obtain the
assertion from the continuity of $F$.\qed

 The Buck's measurable sequences
$\{v_1(n)\},\{v_2(n)\}, \dots, \{v_r(n)\} $ are ca- lled {\it
independent} if and only if for every $x_1,\dots,x_r$ - real numbers
we have
$$
\mu(\bigcap_{j=1}^r \{n\in \N; v_j(n) < x_j\}) =
\prod_{j=1}^r\mu(\{n\in \N; v_j(n) < x_j\}).
$$

\begin{example} We come back to the Example \ref{vdc}. Consider that
the sequences $\{Q_k^{(j)}\}$ are given such that $Q_0^{(j)}=1,
j=1,\dots ,r$ and $Q_k^{(j)}|Q_{k+1}^{(j)}$ for $j=1,\dots ,r$ and
$k=0,1,2 \dots$. If $Q_k^{(j)}, Q_k^{(j_1)}$ are relatively prime
for $j \neq j_1$. Denote $\{\gamma_j(n)\}$ the van der Corput
sequence with base $Q_k^{(j)}$ for $j=1,\dots ,r$. Then these
sequences are independent (see [IPT]).
\end{example}

\begin{theorem}
\label{indep} {\bf Let $\{v_1(n)\},\{v_2(n)\}, \dots \{v_k(n)\}$ be
independent Buck's measurable polyadicly continuous sequences such
that for every $x$ - real number
$$
\mu(\{n\in \N; v_j(n) < x\} =F_j(x)
$$
where $F_j, j=1\dots,k$ are continuous functions defined on real
line. Then the random variables $\tilde{v_1}, \dots, \tilde{v_k}$
are independent.}
\end{theorem}

{\bf Proof.} For $x_1,\dots x_k$ - real numbers we have
$$
\{\alpha \in \Omega; \tilde{v_1}(\alpha) < x_1 \land \dots \land
\tilde{v_k}(\alpha) < x_k\} \subset
$$
$$
\subset cl(\{n\in \N; v_1(n) \le x_1 \land \dots \land v_k(n) \le
x_k\}).
$$
Thus $P(\tilde{v_1}<x_1 \land \dots \land \tilde{v_k}<x_k) \le
F_1(x) \dots F_k(x_k)$ and so from Theorem \ref{thm2} we get
$P(\tilde{v_1}<x \land  \dots \land \tilde{v_k}<_k)\le
P(\tilde{v_1}<x) \dots P(\tilde{v_k}<x_k)$.

From the other hand side we have
$$
 P(\tilde{v_1}\le x_1) \dots P(\tilde{v_k}\le x_k)=
$$
$$
= \mu(\{n\in \N; v_1(n) \le x\}) \dots \mu(\{n\in \N; v_k(n) \le
 x_k\})=
$$
$$
=P(cl(\{n\in \N; v_1(n) \le x_1 \land \dots \land v_k(n) \le x_k\})
\le P(\tilde{v_1}\le x_1 \land  \dots \land \tilde{v_k}\le x_k).
$$
The proof is complete. \qed

\begin{theorem}
\label{CONV} {\bf Let $\{v_1(n)\}, \{v_2(n)\}$ be two independent
Bu- ck's measurable polyadi- cly continuous sequence, such that the
functions
$$
F_1(x)=\mu(\{n\in \N; v_1(n) < x\}), F_2(x)=\mu(\{n\in \N; v_2(n) <
x\})
$$
are continuous. Then $\{v_1(n)+v_2(n)\}$ is Buck's measurable, and}
\begin{equation}
\label{convolution} \mu(\{n\in \N; v_1(n)+v_2(n) < x\}) =
\int_{-\infty}^\infty F_1(x-y)dF_2(y).
\end{equation}
\end{theorem}

{\bf Proof.} Theorem \ref{thm2} implies that the  random variables
$\tilde{v_1}, \tilde{v_2}$ have the distribution functions
$F_1,F_2$. From Theorem \ref{indep} we get that these random
variables are independent and so  the distribution function of
$\tilde{v_1} + \tilde{v_2}$ is given by the integral on right hand
side in (\ref{convolution}). Clearly this function is continuous and
from this we obtain the assertion from Theorem \ref{thm1}. \qed

Thus we get immediately from Theorem \ref{CONV}:

\begin{corollary} {\bf If $\{v_1(n)\}, \{v_2(n)\}$ are two
independent Buck's uniformly distributed sequences then for the
function
$$
G(x) =\mu(\{n\in \N; v_1(n)+v_2(n) <x\})
$$
we have $G(x)=0, x\le 0, G(x)=\frac{x^2}{2}, x \in [0,1], G(x)=2x
-\frac{x^2}{2} - 1, x \in [1,2], G(x)=1, x>2$.}
\end{corollary}

Analogously we can prove:
\begin{corollary} {\bf If $\{v_1(n)\}, \{v_2(n)\}$ are two
independent Buck's uniformly distributed sequences then for the
function
$$
G(x) =\mu(\{n\in \N; v_1(n)v_2(n) <x\})
$$
we have $G(x)=0, x\le 0, G(x)= x(1-\ln x), x \in (0,1), G(x)=1,
x>1$.}
\end{corollary}

{\bf Proof.} Clearly $\tilde{v_1v_2}=\tilde{v_1}\tilde{v_2}$. Then
we have
$$
P(\tilde{v_1v_2}<x) = \int \int _A 1 dF(t_1)dF(t_2),
$$
where $A=\{[t_1,t_2]; t_1t_2 <x; 0\le t_1, t_2 \le 1 \}$.

From Theorem \ref{thm2} we can prove by induction:

\begin{theorem} {\bf If $\{v_1(n)\}, ...,\{v_k(n)\}$ are independent polyadicly
continuous sequences, such that for $j=1,\dots, k$ the functions
$$
F_j(x)=\mu(\{n \in \N; v_j(n) < x\}
$$
are continuous, then $\{v(n)\}=\{v_1(n)+\dots v_k(n)\}$ is
polyadicly continuous function and
$$
F(x)=\mu(\{n \in \N; v(n) < x\}
$$
is continuous.}
\end{theorem}

Let $h: \Omega \to (-\infty, \infty)$ be a continuous function.
Since $\Omega$ is a compact space it is uniformly continuous.
Consider a $m \in \N$. To the function $h$ we can associate a
periodic function $h_m$ with period $m$ in following way:
$$
h_m(\alpha) = h(s) \Longleftrightarrow \alpha \in s+m\Omega.
$$
Clearly
\begin{equation}
\label{intm} \int h_m dP = \frac{1}{m} \sum_{s=0}^{m-1} h(s).
\end{equation}
We have that $\lim_{N \to \infty}\mathfrak(N!,0)=0$ and so from the
uniform continuity of $h$ we get that $h_{N!}$ converges uniformly
to $h$. From (\ref{intm}) we get
\begin{equation}
\label{intinf} \int h dP = \lim_{N \to \infty} \frac{1}{N!}
\sum_{s=0}^{N!-1} h(s).
\end{equation}
The function $h$ restricted on $\N$ is polyadicly continuous. Thus
exists the proper limit $\lim_{m=0}\frac{1}{m}\sum_{s=0}^{m-1}
h(s)$. And so from (\ref{intinf}) we can conclude
\begin{equation}
\label{INT} \int h dP = \lim_{m \to
\infty}\frac{1}{m}\sum_{s=0}^{m-1}h(s).
\end{equation}

\begin{remark}
Let $\{v(n)\}$ be a polyadicly continuous sequence. Then there
exists
\begin{equation}
\label{mean}
 E(v):= \lim_{N \to \infty} \frac{1}{N} \sum_{n=1}^N v(n).
\end{equation}
If the random variable $\tilde{v}$ has continuous distribution
function $F$ then
$$
\int \tilde{v} dP = \int_{-\infty}^{\infty} x dF(x) = E(v).$$
\end{remark}

A sequence of positive integers $\{k_n\}$ is called {\it uniformly
distributed in} $\Z$ if and only for each $m \in \N, r \in \Z$ we
have that $\{n \in \N; k_n \equiv r \mod{m}\} \in \D$ and $d(\{n \in
\N; k_n \equiv r \mod{m}\}) = \frac{1}{m}$. Let us remark that this
type of uniform distribution is firstly defined in [NIV].

In [PAS3] and [PAS4] is proven that for each sequence $\{k_n\}$
uniformly distributed in $\Z$ and polyadicly continuous sequence
$\{v(n)\}$ we have
\begin{equation} \label{INT2}
\lim_{N\to \infty} \sum_{n=1}^N v(k_n)= E(v).
\end{equation}
\begin{theorem} {\bf Let $\{v_1(n)\},\dots,\{v_k(n)\}$ be polyadicly
continuous independent sequences. Then for every functions
$g_1,\dots,g_k$ continuous on real line we have
$$
 \lim_{N \to \infty} \frac{1}{N}
\sum_{n=1}^{N} g_1(v_1(k_n))\dots g_k(v_k(k_n)) = $$ $$= \lim_{N \to
\infty} \frac{1}{N^k}
\Big(\sum_{n=1}^{N}g_1(v_1(k_n))\Big)\dots\Big(\sum_{n=1}^{N}g_k(v_k(k_n))\Big)
$$
for each sequence $\{k_n\}$ uniformly distributed in $\Z$.}
\end{theorem}

{\bf Proof.} If $\{v(n)\}$ is a polyadicly continuous function, then
it is bounded. Every continuous function $g$ defined on real line is
uniformly continuous on closed interval $[b_1, b_2]$ where $b_1$ is
lower bound of  the sequence $\{v(n)\}$ and $b_2$ its upper bound.
Thus the sequence  $\{g(v(n))\}$ is polyadicly continuous also.

Let us consider $\{v_1(n)\},\dots,\{v_k(n)\}$ - polyadicly
continuous independent sequences. Then
$\{g_1(v_1(n))\},\dots,\{g_k(v_k(n))\}$ are polyadicly continuous
and independent also. Thus
$$
E(\{g_1(v_1(n))\}\dots \{g_k(v_k(n))\}) = E(\{g_1(v_1(n))\})\dots
E(\{g_k(v_k(n))\})
$$
and the assertion follows from (\ref{INT2}). \qed

Analogously can be proved that for dispersion $D^2(v)=E((v-E(v))^2)$
the equation
\begin{equation}
\label{dispersion} D^2(v)=D^2(\tilde{v})
\end{equation}
holds. Thus from above the Chebyshev inequality follows:

\begin{theorem} {\bf If $\{v(n)\}$ is polyadicly continuous Buck's measurable
sequence, such that the function
$$
F(x)=\mu(\{n\in \N; v(n) < x\})
$$
is continuous then
$$
\mu(\{n\in \N;|v(n)-E(v)| \ge \varepsilon \})\le
\frac{D^2(v)}{\varepsilon^2}.
$$
for $\varepsilon>0$.}
\end{theorem}

Directly from central limit theorem we get:
\begin{theorem} {\bf Let $\{v_k(n)\}, k=1,2,3 \dots$ be a sequence
of polyadicly continuous sequences that for every $k=1,2,3, \dots$
the sequences $\{v_j(n)\}, j=1\dots k$ are independent and there
exists a continuous function $F$ such that
$$
\mu(\{n\in \N; v_k(n) \le x\}) = F(x), k=1,2,3 \dots
$$
for each real number $x$. Put $E=E(v_k), D^2 = D^2(v_k), k=1,2,3,
\dots$. Then for every $x$ - real number we have
$$
\lim_{k \to \infty}\mu\Big(\Big\{ n \in \N; \frac{v_1(n)+\dots
+v_k(n) - kE}{\sqrt{k}D} \le x \Big\}\Big) =
\frac{1}{\sqrt{2\pi}}\int_{-\infty}^x e^{\frac{-t^2}{2}} dt.
$$}
\end{theorem}


\begin{thebibliography}{TMF9}
 \addcontentsline{toc}{section}{\quad\ \ Bibliography}

\bibitem [BUC]{BUC}
{{\sc Buck, R., C.,}} \textsl{The measure theoretic approach to
density},
 Amer. J. Math \textbf{68},
1946, 560--580

\bibitem [D-T]{D-T}
{{\sc Drmota, M., Tichy, R. F.,}} \textsl{ Sequences, Discrepancies
and Applications, Springer, Berlin Heidelberg},
 Springer, Berlin Heidelberg,
1997

\bibitem [G]{G}
{{\sc Grekos, G.,}} \textsl{On various definitions of density (a
survey)}, Tatra Mt. Math. Publ.,  31, 2005, 17--27
\bibitem [G1]{G1}
{{\sc Grekos, G.,}} \textsl{The density set (a survey)}, Tatra Mt.
Math. Publ.,  31, 2005, 103--111
\bibitem [GST] {GST}
{{\sc Grabner, P. J., Strauch, O., Tichy, R. F.,}}
\textsl{Lp-discrepancyandstatisticalindependenceofsequences, }

Czechoslovak Mathematical Journal,Vol.49(1999),No.1,97–110



\bibitem [IPT]{IPT}
{{\sc  Iaco, M. R., Pasteka, M.,  Tichy R., F.,}} \textsl{Measure
density for set decompositions and uniform distribution} Rend. Circ.
Math. Palermo (2) , 64, No. 2,  2015 , 323 -- 339


\bibitem [K-N]{K-N}
{{\sc  Kuipers, L.,  Niederreiter, H., }} \textsl{ Uniform
distribution of Sequences}, John Wiley and Sons, N.Y. London, Sydney
Toronto, 1974


\bibitem [N]{N}
{{\sc Novoselov, E. V.,}} \textsl{Topological theory of polyadic
numbers}, Trudy Tbilis. Mat. Inst. 27, 1960, 61 -- 69, (in russian)
\bibitem [N1]{N1}
{{\sc Novoselov, E. V.,}} \textsl{New method in the probability
number theory}, Doklady akademii nauk. ser. matem. No. 2, 28, 1964,
307 -- 364,  in russian

\bibitem [NAR]{NAR}
 {{\sc Narkiewicz, W. , }}
 \textsl{Teoria liczb}, (in polish)
 PWN, Warszawa, 1991

\bibitem [NIV]{NIV}
 {{\sc Niven, I.,}}
 \textsl{Uniform distribution of sequences of integers},
 Trans. Amer. Math. Soc. 98, 52 -- 61


\bibitem [PAS]{PAS}
{{\sc  Pa\v st\'eka, M.,}} \textsl{Some properties of Buck's measure
density\index{measure density}},
  Math. Slovaca 42, no. 1,  1992, 15--32

\bibitem [PAS2]{PAS2}
{{\sc  Pa\v st\'eka, M.,}} \textsl{ Remarks on one type of uniform
distribution} Unif. Distrib. Theory 2, No. 1, 2007, 79--92

\bibitem [PAS3]{PAS3}
{{\sc  Pa\v st\'eka, M.,}} \textsl{On four approaches to density}
Spectrum Slovakia 3. Frankfurt am Main: Peter Lang; Bratislava:
VEDA, Publishing House of the Slovak Academy of Sciences, 2014


\bibitem [PAS4]{PAS4}
{{\sc  Pa\v st\'eka, M.,}} \textsl{Remarks on Buck's measure
density\index{measure density}.} Tatra Mt. Math. Publ. 3, 1993,
191--200

\bibitem [P-T]{P-T}
{{\sc Pasteka, M., Tichy, R.,}} \textsl{A note on the correlation
coefficient of arithmetic functions }, Acta Acad. Paed. Agriensis,
Sectio Mathematicae  30, 2003, 109--114


\bibitem [SP]{SP}
{{\sc Strauch, O. Porubsk\'y, \v S.,}} \textsl{  Distribution of
Sequences a Sampler, Peter Lang, SAV, Frankfurt am Main}, Peter
Lang, SAV, Frankfurt am Main, 2005

\bibitem [WEY]{WEY}
 {{\sc Weyl, H.}}
 \textsl{Uber die Gleichverteilung von Zahlen mod. Eins},
 Math. Ann, 77, 1916, 313--352

\end{thebibliography}
\end{document}